\documentstyle[amscd,amssymb,verbatim,amsopn,12pt]{amsart}

\theoremstyle{plain}
\newtheorem{Thm}{Theorem}[section]

\newtheorem{Prop}[Thm]{Proposition}
\newtheorem{Cor}[Thm]{Corollary}
\newtheorem{Con}[Thm]{Conjecture}
\theoremstyle{definition}
\newtheorem{Rem}[Thm]{Remark}

\numberwithin{equation}{section}

\DeclareMathOperator{\Aut}{Aut}

\DeclareMathOperator{\OS}{OS}
\DeclareMathOperator{\ord}{ord}

\newcommand{\bnum}{\begin{enumerate}}
\newcommand{\enum}{\end{enumerate}}
\newcommand{\Z}{\mathbb{Z}}

\begin{document}

\title[POS GROUPS REVISITED]%
	{POS GROUPS REVISITED}

\author[ASHISH KUMAR DAS]%
	{ASHISH KUMAR DAS}

\date{}
\address{Department of Mathematics, North-Eastern Hill
University, Permanent Campus, Shillong-793022, Meghalaya, India.} 
\email{akdas@@nehu.ac.in}


\begin{abstract}
 A finite group $G$ is said to be a POS-group if for each $ x $ in $G$ the cardinality of the set $\{y \in G | o(y) =o(x)\}$ is a divisor of the order of $G$. In this paper we study some of the properties of arbitrary POS-groups, and construct a couple of  new families of nonabelian POS-groups.  We also prove that  the alternating group $A_n$, $n \ge 3$, is not a POS-group. 
\end{abstract}
 
\maketitle

\section{Introduction} \label{S:intro}
  Throughout this paper $G$  denotes a finite group, $o(x)$  the order of a group element $x$, and $|X|$  the cardinality of a set $X$. Also, given a positive integer $n$ and a prime $p$, $\ord_p n$  denotes the largest nonnegative integer $k$ such that $p^k |n$. As in \cite{fJ02},   the order subset (or, order class) of $G$ determined by an element $x \in G$ is defined to be the set $\OS(x) = \{y \in G | o(y) =o(x)\}$. Clearly, $\forall x \in G$, $\OS(x)$ is a disjoint union of some of the conjugacy classes in $G$. The group $G$ is said to have perfect order subsets (in short, $G$ is called a POS-group) if $|\OS(x)|$ is a divisor of $|G|$ for all $x \in G$.  

The object of this paper is to study some of the properties of arbitrary POS-groups, and construct a couple of new families of nonabelian POS-groups. In the process, we re-establish the facts that there are infinitely many nonabelian POS-groups other than the symmetric group $S_3$, and that if a POS-group has its order divisible by an odd prime then it is not necessary that $3$ divides the order of the group (see \cite{fJ02}, \cite{fJ03} and \cite{lT03}). Finally, we prove that  the alternating group $A_n$,  $n \ge 3$, is not a POS-group   (see \cite{fJ03}, Conjecture $5.2$).

\section{Some necessary conditions} \label{S:cond}

Given a positive integer $n$, let $C_n$ denote the cyclic group of order $n$. Then, we have the following characterization for the cyclic POS-groups.
\begin{Prop}\label{pro2.1}
$C_n$ is a POS-group if and only if \,$n=1$ or \,$n = 2^{\alpha} 3^{\beta}$ where  $\alpha \ge 1, \; \beta \ge 0$.
\end{Prop}

\begin{pf}
For each positive divisor $d$ of $n$, $C_n$ has exactly $\phi(d)$ elements of order $d$,  where $\phi$ is the Euler's phi function. So, $C_n$ is a POS-group if and only if $\phi(d)|n$  $\forall$  $d|n$, $i.e.$ if and only if $\phi(n)|n$; noting that $\phi(d)| \phi(n)$ $\forall$ $d|n$.  Elementary calculations reveal that  $\phi(n)|n$ if and only if  \,$n=1$ or \,$n = 2^{\alpha} 3^{\beta}$ where  $\alpha \ge 1, \; \beta \ge 0$. Hence, the proposition follows.
\end{pf}

The following proposition plays a very crucial role in the study of POS-groups (abelian as well as nonabelian).
\begin{Prop}\label{pro2.2}
For each $x \in G$, $|\OS(x)|$  is a multiple of $\phi(o(x))$. 
\end{Prop}

\begin{pf} Define an equivalence relation $\sim$ by setting $a\sim b$ if $a$ and $b$ generate the same cyclic subgroup of $G$. Let $[a]$ denote the equivalence class of $a$ in $G$ under this relation. Then, $\forall x \in G$ and $\forall a \in \OS(x)$, we have  $[a] \subset \OS(x)$, and $|[a]| = \phi(o(a)) = \phi(o(x))$.  Hence it follows that $|\OS(x)| = k\cdot \phi(o(x))$ for all $x \in G$, where $k$ is the number of distinct equivalence classes that constitute $\OS(x)$.
\end{pf}

As an immediate consequence we have the following generalization to the Proposition $1$ and Corollary $1$ of \cite{fJ02}.
\begin{Cor}\label{cor2.3}
If $G$ is a POS-group then, for every prime divisor $p$ of $|G|$, $p-1$ is also a divisor of $|G|$. In particular, every  nontrivial POS-group is of even order.
\end{Cor}

\begin{pf}
By Cauchy's theorem (see \cite{dR96}, page $40$), $G$ has an element of order $p$. So, $G$ being a POS-group, $\phi(p) = p-1$ divides   $|G|$. 
\end{pf}

A celebrated theorem of Frobenius  asserts that if $n$ is a positive divisor of $|G|$ and $X =\{g \in G | g^n = 1\}$, then $n$ divides $|X|$ (see, for example, Theorem 9.1.2 of \cite{mH59}). This result enables us to characterize the $2$-groups having perfect order subsets.
\begin{Prop}\label{pro2.4}
A $2$-group is a POS-group if and only if it is cyclic.
\end{Prop}

\begin{pf}
By Proposition \ref{pro2.1}, every cyclic $2$-group is a POS-group. So, let $G$ be a POS-group with $|G|=2^m$, $m \ge 0$. For $0 \le n \le m$, let $X_n =\{g \in G \, | \, g^{2^n} = 1\}$. Clearly, $X_{n-1} \subset X_n$ for $1 \le n \le m$. We use inductuion to show that $|X_n| = 2^n$ for all $n$ with $0 \le n \le m$. This is equivalent to saying that $G$ is cyclic.  Now, $|x_0|=1=2^0$. So, let us assume that $n \ge 1$. Since $G$ is a POS-group, and since $X_n - X_{n-1}= \{g \in G \, | \, o(g) = 2^n \}$, we have, using Proposition \ref{pro2.2},
\begin{equation}\label{eq2.1}
|X_n| - |X_{n-1}| = |X_n - X_{n-1}| = 0 \text{ or } 2^t
\end{equation}
for some $t$ with $n-1 \le t \le m$. By induction hypothesis, $|X_{n-1}| = 2^{n-1}$, and, by Frobenius' theorem, $2^n$ divides $|X_n|$. Hence, from (\ref{eq2.1}), it follows that $|X_n| = 2^n$. This completes the proof.
\end{pf}

The possible odd prime factors of the order of a nontrivial POS-group are characterized as follows
\begin{Prop}\label{pro2.5}
Let $G$ be a nontrivial POS-group. Then, the odd prime factors (if any) of $|G|$ are of the form $1+2^k t$, where $k \le \ord_2 |G|$ and $t$ is odd, with the smallest one being a Fermat's prime.
\end{Prop}

\begin{pf}
Let $p$ be an odd prime factor of $|G|$. Then, by Corollary \ref{cor2.3}, 
\[
\; p-1 \text{ divides } |G| \; \Longrightarrow  \; \ord_2 (p-1) \le \ord_2 |G|,
\]
which proves the first part. In particular, if $p$ is the smallest odd prime factor of $|G|$ then  $ p-1 = 2^k$,  for some $k \le \ord_2 |G|$. Thus $p = 1 + 2^k$ is a Fermat's prime; noting that $k$ is a power of $2$ as $p$ is a prime. 
\end{pf}

We now determine, through a series of propositions, certain necessary conditions for a group to be a POS-group.
\begin{Prop}\label{pro2.6}
Let $G$ be a nontrivial POS-group with $\ord_2 |G| = \alpha$. If $x \in G$ then the number of distinct odd prime factors in $o(x)$ is at the most $\alpha$. In fact, the bound gets reduces by $(k-1)$ if $\ord_2 o(x) = k \ge 1$.
\end{Prop}

\begin{pf}
If $o(x)$ has $r$ distinct odd prime factors then $2^r |\phi(o(x))$, and so $r \le \alpha$. In addition, if $\ord_2 o(x) = k \ge 1$ then $2^{r+k-1} |\phi(o(x))$, and so $r \le \alpha - (k-1)$.
\end{pf}

\begin{Prop}\label{pro2.7}
If $|G| =2k$ where $k$ is an odd positive integer having at least three distinct prime factors, and if all the Sylow subgroups of $G$ are cyclic, then $G$ is not a POS-group.
\end{Prop}

\begin{pf}
By (\cite{dR96}, $10.1.10$, page $290$ ), $G$ has the following presentation:
\[
G = \langle x,y | x^m = 1 = y^n, x y x^{-1} = y^r \rangle
\]
where $0 \le r < m$, $r^n \equiv 1 \pmod m$, $m$ is odd,  $\gcd(m, n(r-1)) =1$, and $mn = 2k$. Clearly, at least one of $m$ and $n$ is    divisible by two distinct odd primes. So, $o(x)$ or $o(y)$ is divisible by at least two distinct odd primes. The result now follows from Proposition \ref{pro2.6}.
\end{pf}

\begin{Prop}\label{pro2.8}
Let $|G| = {p_1}^{\alpha_1}{p_2}^{\alpha_2}{p_3}^{\alpha_3} \dots {p_k}^{\alpha_k}$ where ${\alpha_1}, {\alpha_2},  \dots, {\alpha_k}$ are positive integers and $2=p_1 <p_2< \dots <p_k$ are primes such that $p_k -1 = {p_1}^{\alpha_1}{p_2}^{\alpha_2}{p_3}^{\alpha_3} \dots {p_{k-1}}^{\alpha_{k-1}}$, $k \ge 2$. If $G$ is a POS-group then the Sylow $p_k$-subgroup of $G$ is cyclic. 
\end{Prop}

\begin{pf}
Note that $G$ has a unique Sylow $p_k$-subgroup, say $P$, so that every element of $G$, of order a power of $p_k$, lies in $P$. Let $m_i$ denote the number of elements of $G$ of order ${p_k}^i$, $1 \le i \le \alpha_k$. Then, by Proposition \ref{pro2.2}, $\phi({p_k}^i)|m_i$. So,
\[
m_i = {p_k}^{i-1} ({p_k} -1)x_i
\]
for some integer $x_i \ge 0$.  If $G$ is a POS-group then we have

\begin{equation}\label{e:1}
x_i | {p_k}^{\alpha_k -i+1} 
\end{equation}
whenever $x_i \neq 0$, $1 \le i \le k$. Now, 
\begin{align}\label{e:2}
&\; \underset{i=1} {\overset{\alpha_k}{\sum}} m_i = |P| - 1 = {p_k}^{\alpha_k} -1 \notag \\
\Rightarrow &\; \underset{i=1} {\overset{\alpha_k}{\sum}} {p_k}^{i-1} \times (x_i -1) = 0  
\end{align}
This gives
\begin{align*}
 &\; x_1 \equiv 1 \pmod {p_k} \\
\Rightarrow &\; x_1 = 1 , \quad \quad \text{by (\ref{e:1}).}
\end{align*}
But, then (\ref{e:2}) becomes
\[
\underset{i=2} {\overset{\alpha_k}{\sum}} {p_k}^{i-1} \times (x_i -1) = 0.
\] 
Repeating the above process inductively, we get
\begin{align*}
&\; x_1 = x_2 = \dots = x_{\alpha_k} = 1 \\
\Rightarrow &\; m_{\alpha_k}= {p_k}^{{\alpha_k} -1} ({p_k} -1) \neq 0.
\end{align*}
This means that  P is cyclic. 
\end{pf}

In view of Proposition \ref{pro2.7} the following corollary is immediate.

\begin{Cor}\label{cor2.9}
If \, $|G| = 42 \times 43^r$, $r \ge 1$, then $G$ is not a POS-group.
\end{Cor}

Finally, we have
\begin{Prop}\label{pro2.10}
Let $G$ be a nontrivial POS-group. Then,  the following assertions hold:
\bnum
\item If $\ord_2 |G| = 1$ then either $|G|=2$, or $3$ divides $|G|$.
\item If $\ord_2 |G| = \ord_3 |G| = 1$ then either $|G|=6$, or $7$ divides $|G|$.
\item If $\ord_2 |G| = \ord_3 |G| = \ord_7 |G| = 1$ then either $|G|=42$, or there exists  a prime $p \ge 77659$ such that $43^2 p$ divides $|G|$.
\enum
\end{Prop}

\begin{pf}
We have already noted that $|G|$ is even. So, let 
\[
|G| =  {p_1}^{\alpha_1} \times {p_2}^{\alpha_2} \times \cdots {p_k}^{\alpha_k}
\]
where $k \ge 1$, $ 2 = p_1<p_2< \dots <p_k$ are primes, and $\alpha_1, \alpha_2, \dots, \alpha_k$ are positive integers.  Now,  for all $i=1, 2, \dots, k$, we have $\gcd(p_i, p_i -1) = 1$. So, in view of Corollary \ref{cor2.3}  we have the following implications:
\begin{align*}
&k \ge 2, \;\alpha_1 = 1 \;\Longrightarrow \; (p_2 -1) | 2 \;\Longrightarrow \; p_2 = 3,\\
&k \ge 3, \;\alpha_1 = \alpha_2 =1 \;\Longrightarrow \; (p_3 -1) | 6 \;\Longrightarrow \; p_3 = 7,\\
&k \ge 4, \;\alpha_1 = \alpha_2 = \alpha_3 = 1 \;\Longrightarrow \; (p_4 -1) | 42 \;\Longrightarrow \; p_4 = 43.
\end{align*}
However,
\[
k \ge 5, \alpha_1 = \alpha_2 = \alpha_3 = \alpha_4 = 1 \;    \Longrightarrow \; (p_5 -1) | 1806 
\]
which is not possible for any prime $p_5 > 43$.  Hence, the theorem follows from Corollary \ref{cor2.9} and the fact that $p=77659$ is the smallest prime greater than $43$ such that $p-1$ divides $2\times 3\times 7\times 43^r$, $r > 1$.
\end{pf}

\begin{Rem}
Using Proposition \ref{pro2.7} and the celebrated theorem of Frobenius one can, in fact, show that if $|G|= 42 \times 43^r \times 77659$,   $r \le 3$, then $G$ is not a POS-group. The proof involves counting of group elements of order powers of $43$.
\end{Rem}

 We have enough evidence in support of  the following conjecture; however, a concrete proof is still eluding.
\begin{Con}
If $G$ is a POS-group such that $\ord_2 |G|$ = $\ord_3 |G|$ = $\ord_7 |G|$ = $1$ then $|G| = 42$.
\end{Con}

\section{Some examples}\label{exam}

Recall (see \cite{dR96}, page $27$) that if $H$ and $K$ are any two groups, and $\theta : H \longrightarrow \Aut(K)$ is a homomorphism then the Cartesian product $H \times K$ forms a group under the binary operation 
\begin{equation}\label{e:4}
(h_1,k_1)(h_2,k_2) = (h_1 h_2, \theta(h_2)(k_1) k_2), 
\end{equation}
where $h_i \in H$, $k_i \in K$,  $i=1,2$. This group is known as the \textit{semidirect product} of $H$ with $K$ (with respect to $\theta$ ), and is denoted by $H \ltimes_\theta K$. Such groups play a very significant role in the construction nonabelian POS-groups.

The following proposition gives a partial characterization of POS-groups whose orders have exactly one distinct odd prime factor.
\begin{Prop}\label{pro3.1}
Let $G$ be a POS-group with $|G|=2^{\alpha} p^{\beta}$ where $\alpha$ and $\beta$ are positive integers, and  $p$ is a Fermat's prime. If $2^{\alpha} < (p-1)^3$ then $G$ is isomorphic to a semidirect product of a group of order $2^{\alpha}$ with the cyclic group $C_{p^\beta}$.
\end{Prop}

\begin{pf}
Since $p$ is a Fermat's prime, $p = 2^{2^k} +1$  where $k \ge 0$.         Let   $X_n =\{g \in G \, | \, g^{p^n} = 1\}$ where $0 \le n \le \beta$. Then, using essentially the same argument as in the proof of Proposition \ref{pro2.4} together with the fact that the order of $2$ modulo $p$ is $2^{k + 1}$, we get $|X_n| = p^n$ for all $n$ with $0 \le n \le \beta$. This implies that $G$ has a unique (hence normal) Sylow $p$-subgroup and it is cyclic. Thus, the proposition follows.
\end{pf}

Taking cue from the above proposition, we now construct a couple new families of nonabelian POS-groups which also serve as counter-examples to the first and the third question posed in section 4 of \cite{fJ02}.
 
\begin{Thm}\label{thm3.2}
Let $p$ be a Fermat's prime. Let $\alpha$, $\beta$ be two positive integers such that $2^\alpha \ge p-1$.  Then there exists a homomorphism $\theta : C_{2^\alpha} \rightarrow \Aut(C_{p^\beta})$ such that the semidirect product $C_{2^\alpha} \ltimes_\theta C_{p^\beta}$ is a nonabelian POS-group.
\end{Thm}

\begin{pf}
Since $p$ is a Fermat's prime, $p = 2^{2^k} +1$  where $k \ge 0$. Also, since the group $U(C_{p^\beta})$ of units in the ring $C_{p^\beta}$ is cyclic and has order $p^{\beta -1} \times 2^{2^k}$, there exists a positive integer $z$ such that
\[
z^{2^{2^k}} \equiv 1 \pmod {p^\beta}, \; \text{and }\;  z^{2^{2^k -1}} \equiv -1 \pmod {p^\beta}.
\]
Moreover, we may choose $z$ in such a way that  
\[
z^{2^{2^k}} \not\equiv 1 \pmod {p^{\beta +1}}
\]
(for, otherwise,  $z$ may be replaced  by $z+ p^\beta$). Let the cyclic groups $C_{2^\alpha}$ and $C_{p^\beta}$ be generated by $a$ and $b$  respectively. Define an automorphism $f : C_{p^\beta} \rightarrow C_{p^\beta}$ by setting $f(b) = b^z$; noting that $\gcd(z,p)=1$. Consider now the homomorphism $\theta : C_{2^\alpha} \rightarrow \Aut(C_{p^\beta})$ defined by $\theta (a) = f.$  
Let $(a^x, b^y) \in C_{2^\alpha} \ltimes C_{p^\beta}$. we may write $x = 2^r m$, $y=p^s n$, where $0 \le r \le \alpha$, $0 \le s \le \beta$, $2 \nmid m$, and $p \nmid n$. It is easy to see that
\begin{equation}\label{e:3}
\theta(a^x)(b^y)= b^{yz^x}.
\end{equation}
So, in $C_{2^\alpha} \ltimes_\theta C_{p^\beta}$, we have, by repeated application of (\ref{e:4}) and (\ref{e:3}),
\begin{equation}\label{e:5}
(a^x, b^y)^{2^{\alpha -r}} = (1, b^\gamma)
\end{equation}
where
\begin{equation}\label{e:6}
\gamma = y \times \dfrac{z^{2^\alpha m} -1}{z^{2^r m} -1}.
\end{equation}
Now, put $c= \ord_p m$.  Then, $m= p^c u$ for some positive integer $u$ such that $p \nmid u$. Therefore, using elementary number theoretic techniques, we have,  for all $r \ge 2^k$, 
\[
z^{2^r m} =  (z^{2^{2^k}})^{2^{r - 2^k} p^c u} \equiv 1 \pmod{p^{\beta +c}} \quad \text{but }\; \not\equiv 1 \pmod{p^{\beta +c +1}}.
\]
On the other hand, if $r < 2^k$  then 
\[
z^{2^r m} \not\equiv 1 \pmod p;
\]
otherwise, since $z$ has order $2^{2^k}$ modulo $p$, we will have 
\[
2^{2^k}|2^r m \; \Longrightarrow  \; 2^k \le r.
\]
Thus, we have
\[
\gamma =
\begin{cases}
p^{\beta +c +s} v, \; &\text{if } r < 2^k, \\
p^s w, \; &\text{if } r \ge 2^k,
\end{cases}
\]
where $v$ and $w$ are two positive integers both coprime to $p$.  This, in turn, means that
\begin{equation}\label{e:7}
o((a^x, b^y)^{2^{\alpha -r}}) =
\begin{cases}
1, \; &\text{if } r < 2^k, \\
p^{\beta -s}, \; &\text{if } r \ge 2^k.
\end{cases}
\end{equation}

Putting $o(a^x, b^y) = t$, we have
\[
(a^x, b^y)^t = (1,1) \;  \Longrightarrow \; a^{tx} = 1 \; 
\Longrightarrow \;  2^\alpha |2^r t m \;  \Longrightarrow \; 2^{\alpha -r} | t, 
\]
since $m$ is odd. Thus, $2^{\alpha -r} | o(a^x, b^y)$. hence, from \ref{e:7}, we have
\begin{equation}\label{e:8}
o(a^x, b^y) =
\begin{cases}
2^{\alpha -r}, \; &\text{if } r < 2^k, \\
2^{\alpha -r} p^{\beta -s}, \; &\text{if } r \ge 2^k.
\end{cases}
\end{equation}

\noindent This enables us to count the number of elements of  $C_{2^\alpha} \ltimes_\theta C_{p^\beta}$   having a  given order, and frame the following table:

\begin{center}
 \begin{tabular}{|c|c|}
\hline Orders of  & Cardinalities of \\
 group elements & corresponding order subsets \\
\hline             $1$               &             $1$     \\
\hline $2^{\alpha -r}$, $(0 \le r < 2^k)$ & $2^{\alpha -r -1} p^\beta$    \\
\hline $2^{\alpha -r}$, $(2^k \le r < \alpha)$ & $ 2^{\alpha -r -1}$ \\     
\hline $p^{\beta -s}$, $(0 \le s < \beta)$ & $p^{\beta -s -1}(p -1)$ \\                                 
\hline $2^{\alpha -r} p^{\beta -s}$,  $ (2^k \le r < \alpha, \; 0 \le s < \beta)$  & $ 2^{\alpha -r -1}  p^{\beta -s -1}  (p -1)$ \\                                               
\hline
\end{tabular}
\end{center}

\noindent It is now easy to see from this table that $C_{2^\alpha} \ltimes_\theta C_{p^\beta}$ is a nonabelian POS-group.  This completes the proof.
\end{pf}

\begin{Rem}
For  $p=5$, taking $z = -1$ in the proof of the above theorem, we get another class of nonabelian POS-groups, namely, $C_{2^\alpha} \ltimes_\theta C_{5^\beta}$ where $\alpha \ge 2$ and $\beta \ge 1$. In this case we have the following table:

\begin{center}
 \begin{tabular}{|c|c|}
\hline Orders of  & Cardinalities of \\
 group elements & corresponding order subsets \\
\hline             $1$               &             $1$     \\
\hline $2^{\alpha}$ & $ 2^{\alpha  -1}  5^\beta$    \\
\hline $2^{\alpha -r}$,  $( 1 \le r < \alpha )$ & $ 2^{\alpha -r -1}$    \\
\hline $5^{\beta -s}$,  $(0 \le s < \beta)$ & $2^2 5^{\beta -s -1}$ \\  
\hline $2^{\alpha -r} 5^{\beta -s}$,  $ (1 \le r < \alpha, \; 0 \le s < \beta)$  & $ 2^{\alpha -r +1}  5^{\beta -s -1}$ \\                         
\hline
\end{tabular}
\end{center}
\end{Rem}

\begin{Rem}\label{r:3}
The argument used in the above theorem also enables us to show that the       semidirect product $C_6 \ltimes_\theta C_7$ is a nonabelian POS-group where the  homomorphism $\theta : C_6 \rightarrow \Aut(C_7)$ is given by $(\theta  (a))(b) = b^2$  (here $a$ and $b$ are generators of $C_6$ and $C_7$ respectively). In this case the element orders are $1,2,3,6,7,14$ and the cardinalities of the corresponding order subsets are $1,1,14,14,6,6$.
\end{Rem}

In (\cite{fJ02}, Theorem 1), it has been proved, in particular, that if ${\Z}_{p^a} \times M$ is a POS-group then  ${\Z}_{p^{a+1}} \times M$ is also a POS-group where $a \ge 1$ and $p$ is a prime such that $p \nmid |M|$. Moreover, as mentioned in the proof of Theorem 1.3 of \cite{fJ03}, the group $M$ need not be abelian. This enables us to construct yet another family of nonabelian POS-groups.
 
\begin{Prop}
Let $M$ be a nonabelian group of order $21$. Then,  $C_{2^a} \times M$ is a POS-group for each $a \ge 1$. 
\end{Prop}

\begin{pf}
In view of the above discussion, it is enough to see that $C_2 \times M$ is a POS-group. In fact,  the element orders and the cardinalities of the corresponding order subsets of $C_2 \times M$ are same as those mentioned in Remark \ref{r:3}
\end{pf}

Finally, we settle Conjecture $5.2$ of \cite{fJ03} regarding $A_n$.

\begin{Prop}
For $n \ge 3$, the alternating group $A_n$ is not a POS-group.
\end{Prop}

\begin{pf}
 It has been proved in \cite{rD72} and also in \cite{jB76} that every positive integer, except $1, 2, 4, 6,$ and $9,$ can be written as the sum of distinct odd primes. Consider now a positive integer $n \ge 3$. It follows that either $n$ or $n-1$ can be written as the sum $p_1 + p_2 + \dots p_k$ where $p_1,  p_2 , \dots, p_k$ are distinct odd primes ($k \ge 1$). Clearly, for such $n$, the number of elements of order $p_1 p_2 \dots p_k$ in $A_n$ is $\dfrac{n!}{p_1 p_2 \dots p_k}$ which does not divide $|A_n| = \dfrac{n!}{2}$. This completes the proof.
\end{pf}

\end{document}